\documentclass{amsart}
\usepackage{amssymb}
\usepackage{mathrsfs}
\usepackage{stmaryrd} 
\usepackage{bbm}
\usepackage{oldgerm}
\usepackage[english]{babel}
\usepackage[T1]{fontenc}
\usepackage[latin1]{inputenc}
\usepackage[all]{xy}
\usepackage{hyperref}
\usepackage{upref}
\usepackage{xcolor}
\usepackage{scrextend}

\hypersetup{
    colorlinks,
    linkcolor={red!50!black},
    citecolor={blue!50!black},
    urlcolor={blue!80!black}
}

\newtheorem{teo}[subsection]{Theorem}
\newtheorem{prop}[subsection]{Proposition}
\newtheorem{cor}[subsection]{Corollary}

\theoremstyle{definition}

\newtheorem{rema}[subsection]{Remark}

\numberwithin{equation}{subsection}

\newcommand{\gtimes}{\stackrel{\leftarrow}{\times}}

\mathchardef\mhyphen="2D

\newcommand{\lgg}{{\ttg}}

\newcommand{\mQ}{{\mathbb Q}}

\newcommand{\mN}{{\mathbb N}}

\newcommand{\mZ}{{\mathbb Z}}

\newcommand{\mG}{{\mathbb G}}

\newcommand{\Et}{{\bf \acute{E}t}}

\newcommand{\bMod}{{\bf Mod}}

\newcommand{\et}{{\rm \acute{e}t}}
\newcommand{\fet}{{\rm f\acute{e}t}}

\newcommand{\rf}{{\rm f}}

\newcommand{\Spec}{{\rm Spec}}

\newcommand{\ob}{{\rm Ob}}

\newcommand{\id}{{\rm id}}

\newcommand{\rE}{{\rm E}}

\newcommand{\rH}{{\rm H}}

\newcommand{\rR}{{\rm R}}

\newcommand{\rW}{{\rm W}}

\newcommand{\oK}{{\overline{K}}}

\newcommand{\oR}{{\overline{R}}}
\newcommand{\oS}{{\overline{S}}}

\newcommand{\oU}{{\overline{U}}}

\newcommand{\oX}{{\overline{X}}}

\newcommand{\ox}{{\overline{x}}}
\newcommand{\oy}{{\overline{y}}}

\newcommand{\oeta}{{\overline{\eta}}}

\newcommand{\ocB}{{\overline{\cB}}}

\newcommand{\uG}{{\underline{G}}}

\newcommand{\uX}{{\underline{X}}}

\newcommand{\ug}{{\underline{g}}}

\newcommand{\uoR}{{\underline{\oR}}}

\newcommand{\uhoR}{{\underline{\hoR}}}

\newcommand{\ulambda}{{\underline{\lambda}}}

\newcommand{\uGamma}{{\underline{\Gamma}}}

\newcommand{\cB}{{\mathscr B}}

\newcommand{\cF}{{\mathscr F}}

\newcommand{\co}{{\mathscr O}}

\newcommand{\cV}{{\mathscr V}}

\newcommand{\fm}{{\mathfrak m}}

\newcommand{\ttg}{{\tt g}}

\newcommand{\hoR}{{\widehat{\oR}}}
\newcommand{\huoR}{\widehat{\uoR}}

\newcommand{\bvg}{{\breve{g}}}

\newcommand{\bvocB}{{\breve{\ocB}}}

\newcommand{\bvpsi}{{\breve{\psi}}}

\newcommand{\bvmZ}{{\breve{\mZ}}}

\newcommand{\tE}{{\widetilde{E}}}

\newcommand{\tG}{{\widetilde{G}}}

\newcommand{\tuG}{{\widetilde{\uG}}}

\begin{document}

\title{The relative Hodge-Tate spectral sequence -- an overview}
\author{Ahmed Abbes and Michel Gros}
\address{A.A. Laboratoire Alexander Grothendieck, CNRS, IHES, Université Paris-Saclay,
35 route de Chartres, 91440 Bures-sur-Yvette, France}
\address{M.G. Univ Rennes, CNRS, IRMAR - UMR 6625, F-35000 Rennes, France}
\email{abbes@ihes.fr}
\email{michel.gros@univ-rennes1.fr}

\begin{abstract}
We give in this note an overview of a recent work \cite{ag} leading to a generalization of the Hodge-Tate spectral sequence to morphisms.
The latter takes place in Faltings topos, but its construction requires the introduction of a relative variant of this topos which is the main novelty of our work. 
\end{abstract}

\maketitle

\setcounter{tocdepth}{1}
\tableofcontents

\section{Introduction}  

\subsection{}\label{intro1}
Let $K$ be a complete discrete valuation field of characteristic $0$, with {\em algebraically closed} residue field of characteristic $p>0$, 
$\co_K$ the valuation ring of $K$, $\oK$ an algebraic closure of $K$, $\co_\oK$ the integral closure of $\co_K$ in $\oK$.
We denote by $G_K$ the Galois group of $\oK$ over $K$, by $\co_C$ the $p$-adic completion of $\co_\oK$, 
by $\fm_C$ the maximal ideal of $\co_C$ and by $C$ its field of fractions.
We set $S=\Spec(\co_K)$ and $\oS=\Spec(\co_\oK)$ and we denote by $s$ (resp.  $\eta$, resp. $\oeta$) 
the closed point of $S$ (resp. generic point of $S$, resp. generic point of $\oS$). 
For any integer $n\geq 0$, we set $S_n=\Spec(\co_K/p^n\co_K)$. For any $S$-scheme $X$, we set
\begin{equation}\label{intro1a}
\oX=X\times_S\oS \ \ \ {\rm and}\ \ \  X_n=X\times_SS_n.
\end{equation}

The following statement, called the {\em Hodge-Tate decomposition}, was conjectured by Tate (\cite{tate} Remark page 180)
and proved independently by Faltings \cite{faltings1,faltings2} and Tsuji \cite{tsuji1,tsuji2}.

\begin{teo}\label{intro2}
For any proper and smooth $\eta$-scheme $X$ and any integer $n\geq 0$, 
there exists a canonical functorial $G_K$-equivariant decomposition
\begin{equation}\label{intro2a}
\rH^n_\et(X_\oeta,\mQ_p)\otimes_{\mQ_p}C\stackrel{\sim}{\rightarrow}\bigoplus_{i=0}^n\rH^i(X,\Omega^{n-i}_{X/\eta})\otimes_KC(i-n).
\end{equation}
\end{teo}

The Hodge-Tate decomposition is equivalent to the existence of a canonical functorial $G_K$-equivariant spectral sequence, the {\em Hodge-Tate spectral sequence},
\begin{equation}\label{intro2b}
\rE_2^{i,j}=\rH^i(X,\Omega^j_{X/\eta})\otimes_KC(-j)\Rightarrow \rH^{i+j}_\et(X_\oK,\mQ_p)\otimes_{\mQ_p}C.
\end{equation}
The two statements are equivalent by a theorem of Tate (\cite{tate} theo.~2).
Indeed, the cohomology group $\rH^0(G_K,C(1))$ vanishes, which implies that the spectral sequence degenerates at $\rE_2$.
The cohomology group $\rH^1(G_K,C(1))$ also vanishes, which implies that the abutment filtration splits. 

This Hodge-Tate spectral sequence, which one can guess implicitly in the work of Faltings \cite{faltings2}, 
has been explicitly formulated only later by Scholze \cite{scholze2}.

\subsection{}\label{intro3}
We give in this note an overview of a recent work \cite{ag} leading to a generalization of the Hodge-Tate spectral sequence to morphisms.
The latter takes place in Faltings topos. Its construction requires the introduction of a relative variant of this topos which is the main novelty of our work. 
Using a different approach, Caraiani and Scholze (\cite{cs} 2.2.4) have constructed a relative Hodge-Tate filtration for proper smooth morphisms of adic spaces. 
Hyodo has also considered earlier a special case for abelian schemes \cite{hyodo2}.  

Beyond the Hodge-Tate spectral sequences, we give in \cite{ag} complete proofs of Faltings' main $p$-adic comparison theorems. 
The latter are essential in the construction of these spectral sequences. 
Although the absolute version of these theorems is rather well-understood, 
the relative version which was very roughly sketched by Faltings in the appendix of \cite{faltings2}, has remained little studied.
Scholze has proved similar results (\cite{scholze1} 1.3 and 5.12) in his setting of adic spaces and pro-étale topos.

In a work in progress, we extend the relative Hodge-Tate spectral sequence to more general coefficients in relation with the 
$p$-adic Simpson correspondence \cite{agt}. 
This sheds new lights on the functoriality of the $p$-adic Simpson correspondence by proper (log)smooth pushforward
(for a related result see the work of Liu and Zhu \cite{lz}). 

\subsection*{Acknowledgments}
We would like first to convey our deep gratitude to G.~Faltings for the continuing inspiration coming from his work on $p$-adic Hodge theory. 
We also thank  very warmly O.~Gabber and T.~Tsuji for the exchanges we had on various aspects discussed in this work. 
Their invaluable expertise has enabled us to avoid long and unnecessary detours.
The first author (A.A) thanks the University of Tokyo and Tsinghua University for their hospitality  
during several visits where parts of this work have been developed and presented. 
He expresses his gratitude to T.~Saito and L.~Fu for their invitations.

\section{The local version of the relative Hodge-Tate spectral sequence}\label{lv}

\subsection{}\label{lv1}
Let $X=\Spec(R)$ be an affine smooth\footnote{We treat in \cite{ag} schemes with toric singularities using logarithmic geometry, 
but for simplicity, we consider in this overview only the smooth case.\label{logsch}} $S$-scheme such that $X_s$ is not empty,
$\oy$ a geometric point of $X_\oeta$. We set $\Gamma=\pi_1(X_\eta,\oy)$ and $\Delta=\pi_1(X_\oeta,\oy)$
and we denote by $(V_i)_{i\in I}$ the universal cover of $X_\oeta$ at $\oy$, or what amounts to the same, 
of its irreducible component containing $\oy$ (\cite{agt} VI.9.7.3).
For each $i\in I$, let $X_i=\Spec(R_i)$ be the normalization of $\oX=X\times_S\oS$ in $V_i$. 
\begin{equation}\label{lv1a}
\xymatrix{
V_i\ar[r]\ar[d]&X_i\ar[d]\\
X_\oeta\ar[r]&\oX}
\end{equation}
The $\co_\oK$-algebras $(R_i)_{i\in I}$ form naturally an inductive system. We denote by $\oR$ its inductive limit, 
\begin{equation}\label{lv1b}
\oR=\underset{\underset{i\in I}{\longrightarrow}}{\lim} \ R_i, 
\end{equation}
and by $\hoR$ its $p$-adic completion, that we equip with the natural actions of $\Gamma$.
The $\Gamma$-representation $\hoR$ is an analogue of the $G_K$-representation $\co_C$. 

Given a proper smooth morphism $g\colon X'\rightarrow X$, 
by analogy with the Hodge-Tate spectral sequence \eqref{intro2b}, one can wonder if there exists a canonical 
$\Gamma$-equivariant spectral sequence 
\begin{equation}\label{lv1c}
\rE_2^{i,j}=\rH^i(X',\Omega^j_{X'/X})\otimes_R\hoR[\frac 1 p](-j)
\Rightarrow \rH^{i+j}_\et(X'_\oy,\mQ_p)\otimes_{\mZ_p}\hoR \  ? 
\end{equation}
We cannot answer this question, but we can construct such a spectral sequence after localizing $X$ at a geometric point of its special fiber.

\subsection{}\label{lv3}
Let $X$ be a smooth\footref{logsch} $S$-scheme, 
$\ox$ a geometric point of $X$ above $s$, $\uX$ the strict localisation of $X$ at $\ox$. 
Let $\oy\rightsquigarrow \ox$ a specialization map, that is an $X$-morphism $u\colon \oy\rightarrow \uX$. 
The latter induces an $X_\oeta$-morphism $\oy\rightarrow \uX_\oeta$. We set $\uGamma=\pi_1(\uX_\eta,\oy)$  
and let $\cV_\ox$ be the category of $\ox$-pointed étale $X$-schemes which are affine.  
For each object $U$ of $\cV_\ox$, we denote by $\oR_U$ the $\co_\oK$-algebra defined as in \eqref{lv1b} for the $S$-scheme $U$ and the geometric point
$\oy\rightarrow U_\oeta$ induced by the canonical morphism $\uX\rightarrow U$. 
The $\co_\oK$-algebras $(\oR_U)_{U\in \cV_\ox}$ form naturally an inductive system. We denote by $\uoR$ its inductive limit,
\begin{equation}\label{lv3a}
\uoR=\underset{\underset{U\in \cV_\ox}{\longrightarrow}}{\lim}\ \oR_U,
\end{equation}
and by $\uhoR$ its $p$-adic completion, that we equip with the natural actions of $\uGamma$.

\begin{teo}[\cite{ag} 6.7.19]\label{lv4}
The assumptions are those of \ref{lv3}, moreover, let $g\colon X'\rightarrow X$ be a smooth projective morphism, $\uX'=X'\times_X\uX$. 
Then, there exists a canonical $\uGamma$-equivariant spectral sequence 
\begin{equation}\label{lv4a}
\rE_2^{i,j}=\rH^i(\uX',\Omega^j_{\uX'/\uX})\otimes_{\co_\uX}\huoR[\frac 1 p](-j)
\Rightarrow \rH^{i+j}_\et(X'_\oy,\mQ_p)\otimes_{\mZ_p}\huoR. 
\end{equation}
\end{teo}

It follows from Faltings' almost purity theorem that the cohomology group
$\rH^0(\uGamma,\huoR\otimes_{\mZ_p}\mQ_p(1))$ vanishes. 
Hence, the spectral sequence \eqref{lv4a} degenerates at $\rE_2$. 
However, the cohomology group $\rH^1(\uGamma,\huoR\otimes_{\mZ_p}\mQ_p(1))$
does not vanish in general. In fact, Hyodo proved already in \cite{hyodo2} that the abutment filtration does not split
in general for abelian schemes.

The spectral sequence \eqref{lv4a} can be globalized over a natural topos whose points are para\-metrized by 
specialization maps $\oy\rightsquigarrow \ox$, from a geometric point $\oy$ of $X_\oeta$ to a geometric point $\ox$ of $X$,
namely {\em Faltings topos}. The latter is at the heart of the Hodge-Tate spectral sequence, even in the absolute case. 
It has been widely studied in (\cite{agt} VI). We briefly review it in the next section.

\begin{rema}
We learned from Scholze that he can answer positively question \eqref{lv1c} 
using the relative Hodge-Tate filtration for proper smooth morphisms of adic spaces that he has developed with Caraiani (\cite{cs} 2.2.4). 
Hyodo has proved the particular case of \ref{lv4} where $X'$ is an abelian scheme over $X$, 
after localizing at a geometric generic point of the special fiber of $X$ \cite{hyodo2}. 
\end{rema}

\section{The global version of the relative Hodge-Tate spectral sequence}\label{ft}

\subsection{}\label{ft1}
Let $X$ be a smooth $S$-scheme\footref{logsch}. 
We denote by $E$ be the category of morphisms $V\rightarrow U$
above the canonical morphism $X_\oeta\rightarrow X$, that is, commutative diagrams 
\begin{equation}\label{ft1a}
\xymatrix{V\ar[r]\ar[d]&U\ar[d]\\
X_\oeta\ar[r]&X}
\end{equation}
such that $U$ is étale over $X$ and the canonical morphism $V\rightarrow U_\oeta$ is {\em finite étale}. 
It is useful to consider the category $E$ as fibered by the functor
\begin{equation}\label{ft1b}
\pi\colon E\rightarrow \Et_{/X}, \ \ \ (V\rightarrow U)\mapsto U,
\end{equation}
over the étale site of $X$. 
The fiber of $\pi$ above an object $U$ of $\Et_{/X}$ is canonically equivalent to the category $\Et_{\rf/U_\oeta}$ of finite étale morphisms 
over $U_\oeta$. We equip it with the étale topology and denote by $U_{\oeta,\fet}$ the associated topos. 
If $U_\oeta$ is connected and if $\oy$ is a geometric point of $U_\oeta$, then the topos $U_{\oeta,\fet}$ is equivalent to
the classifying topos of the profinite group $\pi_1(U_\oeta,\oy)$, {\em i.e.}, the category of discrete sets equipped with a continuous left action
of $\pi_1(U_\oeta,\oy)$.

We equip $E$ with the {\em covanishing} topology (\cite{agt} VI.1.10), that is the topology generated by coverings 
$\{(V_i\rightarrow U_i)\rightarrow (V\rightarrow U)\}_{i\in I}$
of the following two types~:
\begin{itemize}
\item[(v)] $U_i=U$ for all $i\in I$ and $(V_i\rightarrow V)_{i\in I}$ is a covering;
\item[(c)] $(U_i\rightarrow U)_{i\in I}$ is a covering and $V_i=V\times_UU_i$ for all $i\in I$. 
\end{itemize}
The resulting site is called {\em Faltings site} of $X$. We denote by $\tE$ and call {\em Faltings topos} of $X$ the topos of sheaves of sets on $E$. 
It is an analogue of the covanishing topos $X_\et\gtimes_{X_\et}X_{\oeta,\et}$ (\cite{agt} VI.4).

To give a sheaf $F$ on $E$ amounts to give:
\begin{itemize}
\item[(i)] for any object $U$ of $\Et_{/X}$, a sheaf $F_U$ of $U_{\oeta,\fet}$, namely the restriction of $F$ to the fiber
of $\pi$ above $U$;
\item[(ii)] for any morphism $f\colon U'\rightarrow U$ of $\Et_{/X}$, a morphism $\gamma_f\colon F_U\rightarrow f_{\oeta*}(F_{U'})$. 
\end{itemize}
These data should satisfy a cocycle condition for the composition of morphisms and  
a gluing condition for coverings of $\Et_{/X}$ (\cite{agt} VI.5.10). Such a sheaf will be denoted by $\{U\mapsto F_U\}$.

There are three canonical morphisms of topos 
\begin{equation}\label{ft2a}
\xymatrix{
&{X_{\oeta,\et}}\ar[d]_-(0.4){\psi}&\\
{X_\et}&{\tE}\ar[l]_-(0.5){\sigma}\ar[r]^-(0.5)\beta&{X_{\oeta,\fet}}}
\end{equation}
such that 
\begin{eqnarray}
\sigma^*(U)&=&(U_\oeta\rightarrow U)^a, \ \ \ \forall \ U\in \ob(\Et_{/X}),\label{intro4a}\\
\beta^*(V)&=&(V\rightarrow X)^a, \ \ \ \forall \ V\in \ob(\Et_{\rf/X_\oeta}),\label{intro4b}\\
\psi^*(V\rightarrow U)&=&V,\ \ \ \forall \ (V\rightarrow U)\in \ob(E),\label{intro4c}
\end{eqnarray}
where the exponent $^a$ means the associated sheaf. 
The morphisms $\sigma$ and $\beta$ are the analogues of the first and second projections of the covanishing topos $X_\et\gtimes_{X_\et}X_{\oeta,\et}$. 
The morphism $\psi$ is an analogue of the co-nearby cycles morphism (\cite{agt} VI.4.13)

Any specialization map $\oy\rightsquigarrow \ox$ from a geometric point $\oy$ of $X_\oeta$ to a geometric point $\ox$ of $X$, 
determines a point of $\tE$ denoted by $\rho(\oy\rightsquigarrow \ox)$ (\cite{agt} VI.10.18). The collection of these points is conservative (\cite{agt} VI.10.21).

\begin{prop}[\cite{ag} 4.4.2]\label{ft6}
For any locally constant constructible torsion abelian sheaf $F$ of $X_{\oeta,\et}$, we have $\rR^i\psi_*(F)=0$ for any $i\geq 1$.
\end{prop}
This statement is a consequence of the fact that for any geometric point $\ox$ of $X$ over $s$, denoting by $\uX$ the strict
localization of $X$ at $\ox$, $\uX_\oeta$ is a $K(\pi,1)$ scheme (\cite{agt} VI.9.21), {\em i.e.},  if $\oy$ is a geometric point of $\uX_\oeta$, 
for any locally constant constructible torsion abelian sheaf $F$ on $\uX_\oeta$ and any $i\geq 0$, we have an isomorphism 
\begin{equation}
\rH^i(\uX_\oeta,F)\stackrel{\sim}{\rightarrow}\rH^i(\pi_1(\uX_\oeta,\oy),F_\oy). 
\end{equation}
This property was proved by Faltings (\cite{faltings1} Lemma 2.3 page 281), generalizing results of Artin (\cite{sga4} XI). 
It was further generalized by Achinger to the log-smooth case (\cite{achinger} 9.5).

\subsection{}\label{ft3}
For any object $(V\rightarrow U)$ of $E$, we denote by $\oU^V$ the integral closure of $\oU$ in $V$ and we set
\begin{equation}\label{ft3a}
\ocB(V\rightarrow U)=\Gamma(\oU^V,\co_{\oU^V}).
\end{equation} 
The presheaf on $E$ defined above is in fact a sheaf (\cite{agt} III.8.16). We write $\ocB=\{U\mapsto \ocB_U\}$ (cf. \ref{ft1}).  
For any étale $X$-scheme $U$ which is affine, 
the stalk of the sheaf $\ocB_U$ of $U_{\oeta,\fet}$ at a geometric point $\oy$ of $U_\oeta$, 
is the representation $\oR_U$ of $\pi_1(U_\oeta,\oy)$ defined in \eqref{lv1b} for $U$.  

For any specialization map $\oy\rightsquigarrow \ox$, we have 
\begin{equation}\label{ft3b}
\ocB_{\rho(\oy\rightsquigarrow \ox)}=\underset{\underset{U\in \cV_\ox}{\longrightarrow}}{\lim}\ \oR_U,
\end{equation}
where $\cV_\ox$ is the category of $\ox$-pointed étale $X$-schemes $U$ which are affine.

\subsection{}\label{ft7}
For any topos $T$, projective systems of objets of $T$ indexed by the ordered set 
of natural numbers $\mN$, form a topos that we denote by $T^{\mN^\circ}$ (\cite{agt} III.7). 

For any integer $n\geq 0$, we set $\ocB_n=\ocB/p^n\ocB$. 
To take into account $p$-adic topology, we consider the $\co_C$-algebra $\bvocB=(\ocB_n)_{n\geq 1}$
of the topos $\tE^{\mN^\circ}$. 
We work in the category $\bMod_{\mQ}(\bvocB)$ of $\bvocB$-modules up to isogeny (\cite{agt} III.6.1), 
which is a global analogue of the category of $\hoR[\frac 1 p]$-representations of $\Delta$ considered in \ref{lv1}.

\begin{teo}[\cite{ag} 6.7.5]\label{ft4}
Let $g\colon X'\rightarrow X$ be a smooth projective morphism. We denote by 
\begin{equation}
\xymatrix{
{X'^{\mN^\circ}_{\oeta,\et}}\ar[r]^{\bvg_\oeta}&{X^{\mN^\circ}_{\oeta,\et}}\ar[r]^\bvpsi&{\tE^{\mN^\circ}}}
\end{equation}
the morphisms induced by $g_\oeta$ and $\psi$ \eqref{ft2a}, 
and by $\bvmZ_p$ the $\mZ_p$-algebra $(\mZ/p^n\mZ)_{n\geq 1}$ of $X'^{\mN^\circ}_{\oeta,\et}$.  
Then, we have a canonical spectral sequence of $\bvocB_\mQ$-modules
\begin{equation}\label{ft4a}
\rE_2^{i,j}=\sigma^*(\rR^ig_*(\Omega^j_{X'/X}))\otimes_{\sigma^*(\co_X)}\bvocB_\mQ(-j)\Rightarrow \bvpsi_*(\rR^{i+j}\bvg_{\oeta*}(\bvmZ_p))\otimes_{\mZ_p}\bvocB_\mQ.
\end{equation}
\end{teo}

The projectivity condition on $g$ is used in \ref{fmct3} below. 
It should be possible to replace it by the properness of $g$. 

The spectral sequence \eqref{ft4a} is called the {\em relative Hodge-Tate spectral sequence}.
We can easily prove that it is $G_K$-equivariant for the natural $G_K$-equivariant structures on the topos and objects that appear. 
We deduce the following. 

\begin{prop}[\cite{ag} 6.7.13]\label{ft5}
Under the assumptions  of \ref{ft4}, the relative Hodge-Tate spectral sequence \eqref{ft4a} degenerates at $\rE_2$.
\end{prop}

\begin{rema}
Using a different approach, 
Caraiani and Scholze have constructed a relative Hodge-Tate filtration for proper smooth morphisms of adic spaces (\cite{cs} 2.2.4). 
\end{rema}

\section{\texorpdfstring{Faltings' main $p$-adic comparison theorems}{Faltings' main p-adic comparison theorems}}

\subsection{}\label{fmct0}
The assumptions and notation of §~\ref{ft} are in effect in this section.
We denote by $\co_{\oK^\flat}$ the limit of the projective system $(\co_\oK/p\co_\oK)_{\mN}$ 
whose transition morphisms are the iterates of the absolute Frobenius endomorphism of $\co_\oK/p\co_\oK$; 
\begin{equation}
\co_{\oK^\flat}= \underset{\underset{\mN}{\longleftarrow}}{\lim}\ \co_\oK/p\co_\oK.
\end{equation}
It is  a perfect complete non-discrete valuation ring of height $1$ and characteristic $p$. 
We fix a sequence $(p_n)_{n\geq 0}$ of elements of $\co_{\oK}$ such that $p_0=p$ and $p_{n+1}^p=p_n$ for any $n\geq 0$. 
We denote by $\varpi$ the associated element of $\co_{\oK^\flat}$ and we set $\xi=[\varpi]-p$  
in the ring $\rW(\co_{\oK^\flat})$ of $p$-typical Witt vectors of $\co_{\oK^\flat}$.
We have a canonical isomorphism 
\begin{equation}
\co_C(1)\stackrel{\sim}{\rightarrow} p^{\frac{1}{p-1}}\xi \co_C.
\end{equation}

\begin{teo}[\cite{faltings2},\cite{ag} 4.8.13]\label{fmct1}
Assume that $X$ is proper over $S$. Let $i,n$ be integers $\geq 0$, 
$F$ a locally constant constructible sheaf of $(\mZ/p^n\mZ)$-modules of $X_{\oeta,\et}$.
Then, the kernel and cokernel of the canonical morphism 
\begin{equation}\label{fmct1a} 
\rH^i(X_{\oeta,\et},F)\otimes_{\mZ_p}\co_C\rightarrow \rH^i(\tE,\psi_*(F)\otimes_{\mZ_p}\ocB)
\end{equation}
are annihilated by $\fm_C$. 
\end{teo}
We say that morphism \eqref{fmct1a} is an {\em almost isomorphism}. 

This is {\em Faltings' main $p$-adic comparaison theorem} from which he derived all comparaison theorems between $p$-adic étale cohomology 
and other $p$-adic cohomologies. It is also the main ingredient in the construction of the absolute Hodge-Tate spectral sequence \eqref{intro2b}.

We revisit in \cite{ag} Faltings' proof of this important result providing more details. 
It is based on Artin-Schreier exact sequence for the ``perfection'' of the ring $\ocB_1=\ocB/p\ocB$.
One of the main ingredients is a structural statement for almost étale $\varphi$-modules on  $\co_{\oK^\flat}$
verifying certain conditions, including an almost finiteness condition in the sense of Faltings.
In our application to the cohomology of Faltings' topos ringed by the ``perfection'' of $\ocB_1$,
the proof of this last condition results from the combination of three ingredients:
\begin{itemize}
\item[(i)] local calculations of Galois cohomology using Faltings' almost-purity theorem (\cite{faltings2}, \cite{agt} II.8.17); 
\item[(ii)] a fine study of almost finiteness conditions for quasi-coherent sheaves of modules on schemes; 
\item[(iii)] Kiehl's result on the finiteness of cohomology of a proper morphism (\cite{kiehl} 2.9'a) (cf. \cite{egr1} 1.4.7).
\end{itemize}

\subsection{}\label{fmct2}
Next, we explain Faltings' construction of the absolute Hodge-Tate spectral sequence \eqref{intro2b}.  
Assume that $X$ is proper over $S$. 
By \ref{ft6}, for any $i,n\geq 0$, we have a canonical  isomorphism
\begin{equation}\label{fmct2a}
\rH^i(X_{\oeta,\et},\mZ/p^n\mZ)\stackrel{\sim}{\rightarrow}\rH^i(\tE,\psi_*(\mZ/p^n\mZ)).
\end{equation} 
It is not difficult to see that the canonical morphism $\mZ/p^n\mZ\rightarrow \psi_*(\mZ/p^n\mZ)$ is an isomorphism. 
Then, by Faltings' main $p$-adic comparison theorem \ref{fmct1}, we have a canonical morphism 
\begin{equation}\label{fmct2b}  
\rH^i(X_{\oeta,\et},\mZ/p^n\mZ)\otimes_{\mZ_p}\co_C\rightarrow \rH^i(\tE,\ocB_n),
\end{equation}
which is an almost isomorphism. To compute $\rH^i(\tE,\ocB_n)$, we use
the Cartan-Leray spectral sequence for the morphism $\sigma\colon \tE\rightarrow X_{\et}$ \eqref{ft2a},
\begin{equation}\label{fmct2c}
\rE_2^{i,j}=\rH^i(X_\et,\rR^j\sigma_*(\ocB_n))\Rightarrow \rH^{i+j}(\tE,\ocB_n). 
\end{equation}
We deduce the absolute Hodge-Tate spectral sequence \eqref{intro2b} using the following global analogue of Faltings' computation of Galois cohomology.

\begin{teo}[\cite{ag} 6.3.8]\label{fmct3} 
There exists a canonical homomorphism
of graded $\co_{\oX_n}$-algebras of $X_{s,\et}$
\begin{equation}\label{fmct3a}
\wedge (\xi^{-1}\Omega^1_{\oX_n/\oS_n})\rightarrow \oplus_{i\geq 0}\rR^i\sigma_*(\ocB_n),
\end{equation}
where $\xi$ is the element of $\rW(\co_{\oK^\flat})$ defined in \ref{fmct0}, 
whose kernel (resp. cokernel) is annihilated by $p^{\frac{2d}{p-1}}\fm_\oK$ (resp. $p^{\frac{2d+1}{p-1}}\fm_\oK$), where $d=\dim(X/S)$.  
\end{teo}

We prove this result using Kummer theory over the special fiber of Faltings ringed topos $(\tE,\ocB)$.

\subsection{}\label{fmct4} 
Let $g\colon X'\rightarrow X$ be a smooth\footref{logsch} morphism.
We associate to $X'$ objects similar to those associated to $X$ in §~\ref{ft} and we equip them with a prime $^\prime$.  
We have a commutative diagram 
\begin{equation}\label{fmct4a}
\xymatrix{
{X'_{\oeta,\et}}\ar[d]_{g_\oeta}\ar[r]^{\psi'}&{\tE'}\ar[d]^{\Theta}\ar[r]^{\sigma'}&{X'_\et}\ar[d]^g\\
{X_{\oeta,\et}}\ar[r]^\psi&{\tE}\ar[r]^{\sigma}&{X_\et}}
\end{equation}
where $\Theta$ is defined, for any object $(V\rightarrow U)$ of $E$, by
\begin{equation}\label{fmct4b}
\Theta^*(V\rightarrow U)=(V\times_{X}X'\rightarrow U\times_XX')^a,
\end{equation}
where the exponent $^a$ means the associated sheaf. 
We have also a canonical ring homomorphism 
\begin{equation}\label{fmct4c}
\ocB\rightarrow \Theta_*(\ocB').
\end{equation}

\begin{teo}[\cite{faltings2} §~6, \cite{ag} 5.7.4] \label{fmct5} 
Assume that $g\colon X'\rightarrow X$ is projective. Let $i,n$ be integers $\geq 0$, $F'$ a locally constant constructible sheaf of 
$(\mZ/p^n\mZ)$-modules of $X'_{\oeta,\et}$. Then, the canonical morphism 
\begin{equation}\label{fmct5a} 
\psi_*(\rR^ig_{\oeta*}(F'))\otimes_{\mZ_p}\ocB\rightarrow \rR^i\Theta_*(\psi'_*(F')\otimes_{\mZ_p}\ocB')
\end{equation}
is an almost isomorphism. 
\end{teo}

Observe that the sheaves $\rR^ig_{\oeta*}(F)$ ($i\geq 0$) are locally constant constructible on $X_\oeta$ 
by smooth and proper base change theorems. 

Faltings formulated this {\em relative version} of his main $p$-adic comparison theorem 
in \cite{faltings2} and he very roughly sketched a proof in the appendix. 
Some arguments have to be modified and the actual proof in \cite{ag} requires much more work.
It is based on a fine study of the local structure of certain almost-étale $\varphi$-modules 
which is interesting in itself (\cite{ag} 5.5.20). 

The projectivity condition on $g$ is used to prove an almost finiteness result for almost coherent modules. 
We rely on the finiteness results of \cite{sga6} instead of those of \cite{kiehl}. 
It should be possible to replace the projectivity condition on $g$ just by the properness of $g$.

\begin{rema}
Scholze has generalized \ref{fmct1} to rigid varieties following the same strategy (\cite{scholze1} 1.3). 
He has also  proved an analogue of \ref{fmct5} in his setting of adic spaces and pro-étale topos (\cite{scholze1} 5.12).
He deduces it from the absolute case using a base change theorem of Huber. 
\end{rema}

\subsection{}\label{fmct6} 
We keep the assumptions of \ref{fmct5}. Let $n$ be an integer $\geq 0$. 
Since the canonical morphism $\mZ/p^n\mZ\rightarrow \psi'_*(\mZ/p^n\mZ)$ is an isomorphism, 
in order to construct the relative Hodge-Tate spectral sequence \eqref{ft4a}, 
we are led by \ref{fmct5} to compute the cohomology sheaves $\rR^q\Theta_*(\ocB'_n)$ ($q\geq 0$). 
Inspired by the absolute case \eqref{fmct2}, 
the problem is then to find a natural factorization of $\Theta$, to which we can apply the Cartan-Leray spectral sequence. 
Consider the commutative diagram of morphisms of topos
\begin{equation}
\xymatrix{
\tE'\ar[rd]^{\sigma'}\ar[d]_{\tau}&\\
{\tE\times_{X_\et}X'_\et}\ar[r]^-(0.4){\pi}\ar[d]_{\Xi}&{X'_\et}\ar[d]^g\\
{\tE}\ar[r]^-(0.5){\sigma}&{X_\et}}
\end{equation}
We prove that the fiber product of topos $\tE\times_{X_\et}X'_\et$ is in fact a {\em relative Faltings topos},
whose definition was inspired by oriented products of topos, beyond the covanishing topos which inspired our definition of the usual Faltings topos.

\section{Relative Faltings topos}\label{rft}

\subsection{}\label{rft3} 
We keep the assumption and notation of \ref{fmct4}.
We denote by $G$ the category of morphisms $(W\rightarrow U\leftarrow V)$
above the canonical morphisms $X'\rightarrow X\leftarrow X_\oeta$, that is, commutative diagrams 
\begin{equation}\label{rft3a}
\xymatrix{W\ar[r]\ar[d]&U\ar[d]&V\ar[l]\ar[d]\\
X'\ar[r]&X&X_\oeta\ar[l]}
\end{equation}
such that $W$ is étale over $X'$, $U$ is étale over $X$ and the canonical morphism $V\rightarrow U_\oeta$ is {\em finite étale}. 
We equip it with the topology generated by coverings 
\begin{equation}
\{(W_i\rightarrow U_i\leftarrow V_i)\rightarrow (W\rightarrow U \leftarrow V)\}_{i\in I}
\end{equation}
of the following three types~:
\begin{itemize}
\item[(a)] $U_i=U$, $V_i=V$ for all $i\in I$ and $(W_i\rightarrow W)_{i\in I}$ is a covering;
\item[(b)] $W_i=W$, $U_i=U$ for all $i\in I$ and $(V_i\rightarrow V)_{i\in I}$ is a covering;
\item[(c)]  diagrams
\begin{equation}\label{rft3b}
\xymatrix{
W'\ar[r]\ar@{=}[d]&U'\ar[d]\ar@{}|\Box[rd]&V'\ar[l]\ar[d]\\
W\ar[r]&U&V\ar[l]}
\end{equation}
where $U'\rightarrow U$ is any morphism and the right square is Cartesian. 
\end{itemize}
The resulting site is called {\em Faltings relative site} of the morphism $g\colon X'\rightarrow X$. 
We denote by $\tG$ and call {\em Faltings relative topos} of $g$ the topos of sheaves of sets on $G$. 
It is an analogue of the oriented product of topos $X'_\et\gtimes_{X_\et}X_{\oeta,\et}$ (\cite{agt} VI.3).

There are two canonical morphisms 
\begin{equation}\label{rft3c} 
\xymatrix{
{X'_\et}&{\tG}\ar[l]_-(0.4){\pi}\ar[r]^(0.4)\lambda&{X_{\oeta,\fet}}},
\end{equation}
defined by 
\begin{eqnarray}
\pi^*(W)&=&(W\rightarrow X\leftarrow X_\oeta)^a, \ \ \ \forall \ W\in \ob(\Et_{/X'}),\\
\lambda^*(V)&=&(X'\rightarrow X\leftarrow V)^a, \ \ \ \forall \ V\in \ob(\Et_{\rf/X_\oeta}),
\end{eqnarray}
where the exponent $^a$ means the associated sheaf. 
They are the analogues of the first and second projections of the oriented product $X'_\et\gtimes_{X_\et}X_{\oeta,\et}$.

If $X'=X$, $\tG$ is canonically equivalent to Faltings topos $\tE$ \eqref{ft1}. 
Hence, by functoriality of relative Faltings topos, we have a natural factorization of $\Theta\colon \tE'\rightarrow \tE$ into 
\begin{equation}\label{rft3d} 
\tE'\stackrel{\tau}{\longrightarrow}\tG \stackrel{\lgg}{\longrightarrow} \tE. 
\end{equation}
These morphisms fit into the following commutative diagram of morphisms of topos
\begin{equation}
\xymatrix{
&\tE'\ar[d]^{\tau}\ar[r]^{\beta'}\ar[ld]_{\sigma'} &{X'_{\oeta,\fet}}\ar[d]^{g_\oeta}\\
X'_\et\ar[d]_{g}\ar@{}|\Box[rd]&\tG\ar[d]^-(0.5){\lgg}\ar[r]^-(0.4){\lambda}\ar[l]_-(0.4){\pi}&{X_{\oeta,\fet}}\\
X_\et&\tE\ar[ur]_{\beta}\ar[l]_{\sigma}&}
\end{equation}
We prove that {\em the lower left square is Cartesian} \eqref{fmct6}. 

We have a canonical morphism 
\begin{equation}
\varrho\colon X'_\et\gtimes_{X_\et}X_{\oeta,\et}\rightarrow \tG.
\end{equation}
To give a point of $X'_\et\gtimes_{X_\et}X_{\oeta,\et}$ amounts to give a geometric point $\ox'$ of $X'$, 
a geometric point $\oy$ of $X_\oeta$ and a specialization map $\oy \rightsquigarrow g(\ox')$. 
We denote (abusively) such a point by $(\oy \rightsquigarrow \ox')$.  
We prove that the collection of points $\varrho(\oy \rightsquigarrow \ox')$ of $\tG$ is conservative. 

\subsection{}\label{rft5} 
Let $\ox'$ be a geometric point of $X'$, 
$\uX'$ be the strict localization of $X'$ at $\ox'$ and $\uX$  the strict localization of $X$ at $g(\ox')$. 
We denote by $\tuG$ the relative Faltings topos of the morphism
$\uX'\rightarrow \uX$ induced by $g$, by $\ulambda\colon \tuG\rightarrow \uX_{\oeta,\fet}$ the canonical morphism \eqref{rft3c}
and by $\Phi\colon \tuG\rightarrow \tG$ the functoriality morphism. 
There is a canonical section $\theta$ of $\ulambda$,
\begin{equation}\label{rft5a}
\xymatrix{
{\uX_{\oeta,\fet}}\ar[r]^{\theta}\ar[rd]_\id&{\tuG}\ar[d]^{\ulambda}\\
&{\uX_{\oeta,\fet}}}
\end{equation}
We prove that the base change morphism induced by this diagram
\begin{equation}\label{rft5b}
\ulambda_*\rightarrow \theta^*
\end{equation}
is an isomorphism. We set
\begin{equation}\label{rft5d}
\phi_{\ox'}=\theta^*\circ \Phi^*\colon \tG\rightarrow \uX_{\oeta,\fet}.
\end{equation}

If $\oy$ is a geometric point of $\uX_{\oeta}$, we obtain naturally a point $(\oy \rightsquigarrow \ox')$ of  $X'_\et\gtimes_{X_\et}X_{\oeta,\et}$.
Then, for any sheaf $F$ of $\tG$, we have a canonical functorial isomorphism
\begin{equation}\label{rft5e}
F_{\varrho(\oy \rightsquigarrow \ox')}\stackrel{\sim}{\rightarrow} \phi_{\ox'}(F)_{\oy}.
\end{equation}

\begin{prop}[\cite{ag} 3.4.32] \label{rft6} 
Under the assumptions of \ref{rft5}, for any abelian sheaf $F$ of $\tG$ and any $q\geq 0$, we have a canonical isomorphism
\begin{equation}
\rR^q\pi_*(F)_{\ox'}\stackrel{\sim}{\rightarrow}\rH^q(\uX_{\oeta,\fet},\phi_{\ox'}(F)).
\end{equation}
\end{prop}

\begin{cor}[\cite{ag} 6.5.17]
Let $(\oy\rightsquigarrow \ox')$ be a point of $X'_\et\gtimes_{X_\et}X_{\oeta,\et}$, 
$\uX'$ the strict localization of $X'$ at $\ox'$, $\uX$ the strict localization of $X$ at $g(\ox')$,
$\ug\colon \uX'\rightarrow \uX$ the morphism induced by $g$, 
\begin{equation}
\varphi'_{\ox'}\colon \tE'\rightarrow \uX'_{\oeta,\fet}
\end{equation}
the canonical morphism analogue of \eqref{rft5d}. 
Then, for any abelian group $F$ of $\tE'$ and any $q\geq 0$, we have a canonical functorial isomorphism
\begin{equation}
(\rR^q\tau_*(F))_{\varrho(\oy\rightsquigarrow \ox')}\stackrel{\sim}{\rightarrow}
\rR^q\ug_{\oeta,\fet*}(\varphi'_{\ox'}(F))_{\oy}.
\end{equation}
\end{cor}

\subsection{}\label{rft7} 
We consider the following ring of $\tG$,
\begin{equation}
\ocB^!=\tau_*(\ocB').
\end{equation}
We have canonical homomorphisms $\ocB\rightarrow \lgg_*(\ocB^!)$ and $\hbar'_*(\co_{\oX'})\rightarrow \pi_*(\ocB^!)$,
where $\hbar'\colon \oX'\rightarrow X'$ is the canonical projection. Hence, we may consider $\lgg$ and $\pi$ as morphisms of ringed topos. 

For any point $(\oy \rightsquigarrow \ox')$ of  $X'_\et\gtimes_{X_\et}X_{\oeta,\et}$, we prove that 
the ring $\ocB^!_{\varrho(\oy \rightsquigarrow \ox')}$ is normal and strictly henselian. Moreover, the canonical homomorphism
$\co_{\oX',\ox'}\rightarrow \ocB^!_{\varrho(\oy \rightsquigarrow \ox')}$ is local and injective. 

\subsection{}\label{rft11} 
Assume that $X=\Spec(R)$ and $X'=\Spec(R')$ are affine. Let $\oy'$ be a geometric point of $X'_\oeta$, $\Delta'=\pi_1(X'_\oeta,\oy')$,
$(W_j)_{j\in J}$ the universal cover of $X'_\oeta$ at $\oy'$, $\oy=g_\oeta(\oy')$, $\Delta=\pi_1(X_\oeta,\oy)$ and
$(V_i)_{i\in I}$ the universal cover of $X_\oeta$ at $\oy$.  
For every $i\in I$, $(V_i\rightarrow X)$ is naturally an object of $E$ and for every $j\in J$, 
$(W_j\rightarrow X')$ is naturally an object of $E'$. We set  
\begin{eqnarray}
\oR&=&\underset{\underset{i\in I}{\longrightarrow}}{\lim}\ \ocB(V_i\rightarrow X),\label{rft11a} \\
\oR'&=&\underset{\underset{j\in J}{\longrightarrow}}{\lim}\ \ocB'(W_j\rightarrow X'). \label{rft11b} 
\end{eqnarray}
We recover the $\co_\oK$-algebras defined in \eqref{lv1b}. We equip them with the natural actions of $\Delta$ and $\Delta'$. 
For  every $i\in I$, there exists a canonical $X'$-morphism $\oy'\rightarrow X'\times_XV_i$. We denote by 
$V'_i$ the irreducible component of $X'\times_XV_i$ containing $\oy'$ and by $\Pi_i$ the corresponding subgroup of $\Delta'$. 
Then, $(V'_j\rightarrow X')$ is naturally an object of $E'$. 
We set 
\begin{eqnarray}
\oR^!&=&\underset{\underset{i\in I}{\longrightarrow}}{\lim}\ \ocB'(V'_i\rightarrow X'),\label{rft11d} \\
\Pi&=&\bigcap_{i\in I}\Pi_i.\label{rft11c}
\end{eqnarray}
We have canonical homomorphisms $\oR\rightarrow \oR^! \rightarrow \oR'$.

For any geometric point $\ox'$ and any specialization map $\oy\rightsquigarrow g(\ox')$, 
we prove that we have a canonical isomorphism (determined by the choice of $\oy'$)
\begin{equation}\label{rft11e} 
\ocB^!_{\varrho(\oy\rightsquigarrow \ox')}\stackrel{\sim}{\rightarrow}
\underset{\underset{\ox'\rightarrow U'\rightarrow U}{\longrightarrow}}{\lim}\ \oR^!_{U'\rightarrow U},
\end{equation}
where the inductive limit is taken over the category of morphisms $\ox'\rightarrow U'\rightarrow U$ over $\ox'\rightarrow X'\rightarrow X$, 
with $U'$ affine étale over $X'$ and $U$ affine étale over $X$, and $\oR^!_{U'\rightarrow U}$ is the corresponding ring \eqref{rft11d}.

\begin{prop}[\cite{ag} 5.2.29]\label{rft12} 
We keep the assumptions of \ref{rft11} and we assume moreover that $g$ fits into a commutative diagram\footref{logsch} 
\begin{equation}
\xymatrix{
X'\ar[r]^-(0.5){\iota'}\ar[d]_g&{\mG_{m,S}^{d'}}\ar[d]^\gamma\\
X\ar[r]^-(0.5){\iota}&{\mG_{m,S}^{d}}}
\end{equation}
where the morphisms $\iota$ and $\iota'$ are étale, $d$ and $d'$ are integers $\geq 0$ and $\gamma$ 
is a smooth homomorphism of tori over $S$. 
Let $n$ be an integer $\geq 0$. Then,
\begin{itemize}
\item[{\rm (i)}] There exists a canonical homomorphism of graded $\oR^!$-algebras
\begin{equation}
\wedge (\Omega^1_{R'/R}\otimes_{R'}(\oR^!/p^n\oR^!)(-1))\rightarrow \oplus_{i\geq 0}\rH^i(\Pi,\oR'/p^n\oR'),
\end{equation}
which is almost injective and its cokernel is killed by $p^{\frac{1}{p-1}}\fm_\oK$. 
\item[{\rm (ii)}] The $\oR^!$-module $\rH^i(\Pi,\oR'/p^n\oR')$ is almost finitely presented for all $i\geq 0$, and it almost vanishes 
for all $i\geq r+1$, where $r=\dim(X'/X)$.
\end{itemize}
\end{prop}

This is a relative version of Faltings' computation of Galois cohomology of $\oR$, that relies on his almost purity theorem 
(\cite{faltings2} Theorem 4 page 192, \cite{agt} II.6.16). 
The statement can be globalized using Kummer theory on the special fiber of the ringed topos $(\tE',\ocB')$ into the following. 

\begin{teo}[\cite{ag} 6.6.4]\label{rft8} 
For any integer $n\geq 1$, there exists a canonical homomorphism of graded $\ocB^!$-algebras of $\tG$
\begin{equation}
\wedge (\pi^*(\xi^{-1}\Omega^1_{\oX'_n/\oX_n}))\rightarrow \oplus_{i\geq 0}\rR^i\tau_*(\ocB'_n),
\end{equation}
where $\pi^*$ denotes the pullback by the morphism of ringed topos $\pi\colon (\tG,\ocB^!)\rightarrow (X'_\et,\hbar'_*(\co_{\oX'}))$, 
whose kernel (resp. cokernel) is annihilated by $p^{\frac{2r}{p-1}}\fm_\oK$ (resp. $p^{\frac{2r+1}{p-1}}\fm_\oK$), where $r=\dim(X'/X)$.  
\end{teo} 
 
\subsection{}\label{rft13} 
Next we consider the Cartan-Leray spectral sequence 
\begin{equation}\label{rft13a}
\rE_2^{i,j}=\rR^i\lgg_*(\rR^j\tau_*(\ocB'_n))\Rightarrow \rR^{i+j}\Theta_*(\ocB'_n).
\end{equation} 
Taking into account \ref{rft8}, to obtain the relative Hodge-Tate spectral sequence \eqref{ft4a}, we need to prove a base change theorem relatively
to the Cartesian diagram 
\begin{equation}\label{rft13b}
\xymatrix{
\tG\ar[r]^{\pi}\ar[d]_{\lgg}&{X'_\et}\ar[d]^g\\
{\tE}\ar[r]^{\sigma}&{X_\et}}
\end{equation}

\begin{teo}[\cite{ag} 6.5.5] \label{rft14} 
Assume $g$ proper. Then, for any torsion abelian sheaf $F$ of $X'_\et$ and any $q\geq 0$, the base change morphism 
\begin{equation}
\sigma^*(\rR^qg_*(F))\rightarrow \rR^q\lgg_*(\pi^*(F))
\end{equation}
is an isomorphism.
\end{teo}

The proof is inspired by a base change theorem for oriented products due to Gabber. It reduces to proper base change theorem 
for étale topos. 

\begin{prop}[\cite{ag} 6.5.29] \label{rft15} 
For any integer $n\geq 0$, the canonical homomorphism
\begin{equation}
\ocB_n\boxtimes_{\co_X}\co_{X'}\rightarrow \ocB^!_n,
\end{equation}
where the exterior tensor product of rings is relative to the Cartesian diagram \eqref{rft13b},
is an almost isomorphism.
\end{prop}

\begin{teo}[\cite{ag} 6.5.30] \label{rft16} 
Assume that the morphism $g$ is proper.
Then, there exists an integer $N\geq 0$ such that for any integers $n\geq 1$ and $q\geq 0$, 
and any coherent $\co_{X'_n}$-module $\cF$ which is $X_n$-flat \eqref{intro1a}, the kernel and cokernel of the base change morphism 
\begin{equation}
\sigma^*(\rR^qg_*(\cF))\rightarrow \rR^q\lgg_*(\pi^*(\cF)),
\end{equation}
where $\sigma^*$ and $\pi^*$ denote the pullbacks in the sense of ringed topos,  are annihilated by $p^N$.
\end{teo}

\begin{prop}[\cite{ag} 6.5.31] \label{rft17} 
Let $n, q$ be integers $\geq 0$, $\cF$ a coherent $\co_{X'_n}$-module which is $X_n$-flat \eqref{intro1a}. 
Assume that the morphism $g$ is proper and that for any integer $i\geq 0$, 
the $\co_{X_n}$-module $\rR^ig_*(\cF)$ is locally free (of finite type). 
Then, the base change morphism
\begin{equation}
\sigma^*(\rR^q g_*(\cF))\rightarrow \rR^q\lgg_*(\pi^*(\cF)),
\end{equation}
where $\sigma^*$ and $\pi^*$ denote the pullbacks in the sense of ringed topos, is an  almost isomorphism.
\end{prop}

\subsection{}
Let $n,q$ be integers $\geq 0$. 
Since the canonical morphism $\mZ/p^n\mZ\rightarrow \psi'_*(\mZ/p^n\mZ)$ is an isomorphism,  
we deduce from \ref{fmct5}, for any $q\geq 0$, a canonical morphism
\begin{equation}
\psi_*(\rR^qg_{\oeta*}(\mZ/p^n\mZ))\otimes_{\mZ_p}\ocB\rightarrow \rR^q\Theta_*(\ocB'_n),
\end{equation}
which is an almost isomorphism. 
The relative Hodge-Tate spectral sequence \eqref{ft4a} is then deduced from the Cartan-Leray spectral sequence 
\eqref{rft13a} using \ref{rft8} and \ref{rft16}.

\end{document}